\documentclass[12pt, a4paper]{article}

\usepackage[utf8]{inputenc}
\usepackage{fixltx2e}

\usepackage{geometry}
 \usepackage{fancyhdr}
\fancypagestyle{customheader}{
  \fancyhf{}
  \fancyhead[R]{%
  \raisebox{3pt}{\includegraphics[height=30pt]{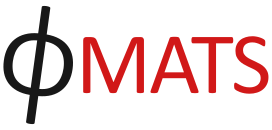}} 
}
  \fancyhead[L]{%
  \vspace{8pt}%
  \begin{minipage}[b]{0.4\textwidth}
    \raggedright
    \small\color{black!80}\textbf{Finite Element Theory for PHIMATS}\\[-0.3em]
    {\color{crimson}\rule{\linewidth}{1pt}}
  \end{minipage}
}
  \fancyfoot[C]{\thepage}
}

\usepackage{afterpage}

\usepackage{graphicx}
\usepackage{amsmath}
\DeclareMathOperator{\Tr}{Tr}
\usepackage{fancybox}
\usepackage[bottom]{footmisc}

\usepackage{algorithm2e,  algpseudocode, amssymb}

\usepackage[
backend=biber,
style=numeric-comp,
sorting=none,
]{biblatex}

\AtBeginDocument{%
  \renewbibmacro{in:}{}%
}

\usepackage{hyperref}
\hypersetup{
    colorlinks=true,
    linkcolor=black,
    citecolor=cyan,
    urlcolor=cyan,
}
\usepackage{xcolor}
\definecolor{crimson}{RGB}{204,  0,  15}
\usepackage{sectsty,titling}
\sectionfont{\color{crimson}}
\subsectionfont{\color{crimson}}
\subsubsectionfont{\color{crimson} \normalfont\itshape}

\title{\textcolor{crimson}{\textbf{Finite Element Theory for PHIMATS}}}

\date{October 2025}

\author{Abdelrahman Hussein \\ email \href{mailto:a.h.a.hussein@outlook.com}{a.h.a.hussein@outlook.com}, \href{mailto:abdelrahman.hussein@oulu.fi}{abdelrahman.hussein@oulu.fi}} 

\definecolor{sub}{HTML}{cde4ff}     
\definecolor{main}{HTML}{5989cf}    
\usepackage[many]{tcolorbox}    	

\newcounter{myBoxCounter}[section]
\renewcommand{\themyBoxCounter}{\thesection.\arabic{myBoxCounter}}

\newtcolorbox[auto counter, number within=section]{boxH}[2][]{
    title=Box~\themyBoxCounter: #2,
    colback = sub, 
    colframe = main, 
    boxrule = 0pt, 
    before skip=20pt plus 2pt,
    after skip=20pt plus 2pt,
    #1,
}

\newcommand{\nextboxnumber}{\stepcounter{myBoxCounter}}

\numberwithin{equation}{section}

\usepackage{bm}

\begin{document}

\maketitle

\noindent\rule{\textwidth}{1pt}

\vspace{1cm}

This document summarizes the main concepts of the finite element (FE) theory and constitutive relations as implemented in the open-source code \textit{phase-field multiphysics materials simulator PHIMATS} (\href{https://github.com/ahcomat/PHIMATS.git}{GitHub Repository}). PHIMATS is written in \texttt{C++} and uses Python for pre- and post-processing. It provides tools for discretizing the weak form of partial differential equations (PDE), interfacing with \texttt{PETSc} \cite{petsc-web-page} data structures (\texttt{Vec}/\texttt{Mat}) and solvers (\texttt{KSP}/\texttt{SNES}). The framework supports both single-physics and coupled multiphysics problems primarily using staggered coupling schemes. Hands-on examples can be found in the \texttt{CaseStudies} directory on GitHub repository. Rather than detailing the derivations of specific models, this document focuses on the key mathematical formulations and numerical strategies used within the implementation. For in-depth theoretical discussions, the reader is encouraged to consult the references. For citing this document, please use: [Abdelrahman Hussein. Finite Element Theory for PHIMATS. 2025. doi: 10.48550/ARXIV.2502.16283]. 

\vspace{0.75cm}

\noindent \rule{17cm}{1pt}

\tableofcontents

\clearpage
\pagestyle{customheader}

\pagebreak

\section{Isoparametric formulation}

The geometry of a domain can be discretized by building blocks called elements, that are defined by the nodes forming these elements and their coordinates. The unknown variables of interest are calculated at these nodes, i.e. nodal unknowns or \emph{degrees of freedom (DOFs)}. These unknowns can be interpolated at any location within an element using the so called (interpolation) shape functions. Similarly, shape functions can be used to interpolate the geometry of the element. For the \emph{isoparametric formulation}, the same shape functions are used to interpolate both the geometry and the nodal unknowns \cite{Biner2017}. 

The shape functions $N_i(\boldsymbol{\xi})$ are constructed to have the value of 1 for a given node and zero everywhere else. $\boldsymbol{\xi}$ is a non-dimensional coordinate system called \emph{local} or \emph{element} coordinate system, which vary in the range [-1,1]. We will use a 4-node quad element as an example for $N_i(\boldsymbol{\xi})$. For a degree of freedom $f$ with known values at the element nodes $f_1, f_2, f_3, f_4$, the value at any point $\boldsymbol{\xi}$ within the element are given by 

\begin{equation}
    f = \sum_i^n N_i(\boldsymbol{\xi}) f_i = N_1(\boldsymbol{\xi}) f_1 + N_2(\boldsymbol{\xi}) f_2 + N_3(\boldsymbol{\xi}) f_3 + N_4(\boldsymbol{\xi}) f_4
\end{equation}

\noindent Where $n$ is the number of nodes\footnote{$N_i=\mathbf{N}$ will be represented as a row vector. And since all the shape functions $N_i(\boldsymbol{\xi})$ have to be evaluated for a given Gauss point $\boldsymbol{\xi}$, they will be stored for each element as $\texttt{nGauss} \times \texttt{nElNodes}$.} per element.

\nextboxnumber
\begin{boxH}[label=ShapeFunctions]{Example of shape functions for a 4-node quadrilateral element}
    \begin{equation*}
        \begin{split}
        N_1 & = 0.25(1 - \xi)(1 - \eta) \\
        N_2 & = 0.25(1 + \xi)(1 - \eta) \\
        N_3 & = 0.25(1 + \xi)(1 + \eta) \\
        N_4 & = 0.25(1 - \xi)(1 + \eta) \\
        \end{split}
    \end{equation*}
\end{boxH}

\noindent Similarly, any point $\boldsymbol{\xi}$ within the element can be mapped to \emph{cartesian} coordinates by 

\begin{equation}
    \begin{split}
    \mathbf{x(\xi)} &= \sum_i^n N_i(\boldsymbol{\xi}) \mathbf{x}_i
    \end{split}
    \label{Eq:CoordsMapping}
\end{equation}

\noindent Where $\mathbf{x}_i$ are the cartesian coordinates of the nodes. 

\subsection{Derivatives}

The cartesian derivatives\footnote{Note this is a matrix of dimensions $\texttt{nDim} \times \texttt{nElNodes}$.} of the variable $f$ at any point $\boldsymbol{\xi}$ defined at the element nodes can be evaluated as

\begin{equation}
    \begin{split}
    \frac{\partial f(\boldsymbol{\xi})}{\partial \mathbf{x}} &= \sum_i^n \frac{\partial N_i(\boldsymbol{\xi})}{\partial \mathbf{x}}f_i^\top \\
    &= \left[
        \begin{array}{rrrr}
        \frac{\partial N_1(\boldsymbol{\xi})}{\partial x} & \dots & \frac{\partial N_n(\boldsymbol{\xi})}{\partial x} \\
        \frac{\partial N_1(\boldsymbol{\xi})}{\partial y} & \dots & \frac{\partial N_n(\boldsymbol{\xi})}{\partial y} \\
        \frac{\partial N_1(\boldsymbol{\xi})}{\partial z} & \dots & \frac{\partial N_n(\boldsymbol{\xi})}{\partial z} \\
        \end{array}
        \right] \left[
        \begin{array}{r}
        f_1 \\
        \vdots \\
        f_n \\
        \end{array}
        \right] \\
    \end{split}
    \label{Eq:IsoparamDeriv}
\end{equation}

\noindent Using the chain rule

\begin{equation}
    \sum_i^n \frac{\partial N_i(\boldsymbol{\xi})}{\partial \mathbf{x}} = \sum_i^n \frac{\partial N_i(\boldsymbol{\xi})}{\partial \boldsymbol{\xi}} \frac{\partial \boldsymbol{\xi}}{\partial \mathbf{x}}
    \label{Eq:ShapeDeriv}
\end{equation}

\noindent $\frac{\partial \mathbf{x}(\boldsymbol{\xi})}{\partial \boldsymbol{\xi}} = \frac{\partial x_i(\boldsymbol{\xi})}{\partial \xi_j}$ is called the jacobian matrix $\mathbf{J}(\boldsymbol{\xi})$ with components

\begin{equation}
    \mathbf{J}(\boldsymbol{\xi}) = \left[\begin{array}{rrr}
        \frac{\partial x}{\partial \xi} & \frac{\partial x}{\partial \eta} & \frac{\partial x}{\partial \zeta} \\
        \frac{\partial y}{\partial \xi} & \frac{\partial y}{\partial \eta} & \frac{\partial y}{\partial \zeta} \\
        \frac{\partial z}{\partial \xi} & \frac{\partial z}{\partial \eta} & \frac{\partial z}{\partial \zeta} \\
        \end{array}\right]
\end{equation}

\noindent Using Eq.~(\ref{Eq:CoordsMapping})

\begin{equation}
    \begin{split}
    \mathbf{J}(\boldsymbol{\xi}) &= \sum_i^n \frac{\partial N_i(\boldsymbol{\xi})}{\partial \boldsymbol{\xi}} \mathbf{x_i}^\top \\
    &= \left[
        \begin{array}{rrrr}
        \frac{\partial N_1(\boldsymbol{\xi})}{\partial x} & \dots & \frac{\partial N_n(\boldsymbol{\xi})}{\partial x} \\
        \frac{\partial N_1(\boldsymbol{\xi})}{\partial y} & \dots & \frac{\partial N_n(\boldsymbol{\xi})}{\partial y} \\
        \frac{\partial N_1(\boldsymbol{\xi})}{\partial z} & \dots & \frac{\partial N_n(\boldsymbol{\xi})}{\partial z} \\
        \end{array}
        \right] \left[
        \begin{array}{rrr}
        x_1 & y_1 & z_1 \\
        \vdots \\
        x_n & y_n & z_n \\
        \end{array}\right]
    \end{split}
\end{equation}

\noindent Substituting in Eq.~(\ref{Eq:ShapeDeriv})

\begin{equation}
    \sum_i^n\frac{\partial N_i(\boldsymbol{\xi})
    }{\partial \mathbf{x}} = J(\boldsymbol{\xi})^{-1} \sum_i^n\frac{\partial N_i(\boldsymbol{\xi})}{\partial \boldsymbol{\xi}} 
\end{equation}

\noindent It is worth noting that the volume/area mapping can be obtained using 

\begin{equation}
    d\mathbf{x} = \det[\mathbf{J(\boldsymbol{\xi})}] d\boldsymbol{\xi} 
\end{equation}

\subsection{Numerical integration}

\noindent Gauss quadrature is used for integration over the element domain (in local coordinates) of the form

\begin{equation}
    \begin{split}
    \int_{-1}^{+1}\int_{-1}^{+1}\int_{-1}^{+1}f(\boldsymbol{\xi})d\boldsymbol{\xi} = \int_{\Omega_e} f(\boldsymbol{\xi})d\boldsymbol{\xi} &\approx \sum^M_{I=1}\sum^M_{J=1}\sum^M_{K=1} w_Iw_Jw_K f(\boldsymbol{\xi}) \\
    &\approx \sum^M_{I=1}\sum^M_{J=1}\sum^M_{K=1} w_Iw_Jw_K f(\xi_I,\eta_J,\zeta_K) \\
    \end{split}
\end{equation} 

\noindent Where $w_{I,J,K}$ is the weight of the the integration point $\boldsymbol{\xi}$. Therefore, the volume integral over all the element can be taken as the summation of volume integral of all its gauss points.

\section{The weak form}
\label{Sec: WeakForm}

The FEM method for solving a PDE starts with forming a variational principle transforming the \emph{strong form} to \emph{weighted-residual} or a \emph{weak form} \cite{Bower2009}. This weak form transforms the PDE to a minimization problem. This involves the following steps \cite{Biner2017, moose-fem-principles}:

\begin{enumerate}
    \item Multiply the PDE with a \emph{test} function. In the Galerkin FEM, the test function is the "trial" function. 
    \item Integrate over the whole domain $\Omega$.
    \item Integrate by parts the terms involving spatial gradients using the divergence theorem.
    \item Provide the given boundary conditions.
  \end{enumerate}

\noindent A general case for integration by parts

\begin{equation*}
    \int_\Omega \psi (\nabla \cdot \nabla u) \, d\Omega = -\int_\Omega \nabla \psi \cdot \nabla u \, d\Omega + \int_{\Gamma} \psi (\nabla u \cdot \mathbf{n}) \, d\Gamma
\end{equation*}

\noindent For example, consider the PDE

\begin{equation}
    \frac{\partial C}{\partial t} + \nabla \cdot D(\mathbf{x}) \, \nabla C = 0
    \label{Eq:Example}
\end{equation}

\noindent Step(1): Multiply Eq.~(\ref{Eq:Example}) by a test function $\psi$

\begin{equation}
    \psi (\frac{\partial C}{\partial t}) + \psi (\nabla \cdot D(\mathbf{x}) \, \nabla C) = 0
\end{equation}

\noindent Step(2): Integrate over the domain 

\begin{equation}
    \int_\Omega \psi (\frac{\partial C}{\partial t})d\Omega + \int_\Omega \psi (\nabla \cdot D(\mathbf{x}) \, \nabla C) \, d\Omega = 0
\end{equation}

\noindent Step(3): Integrate by parts the terms involving gradients, i.e. the second term  of the LHS

\begin{equation}
    \int_\Omega \psi (\nabla \cdot D(\mathbf{x}) \, \nabla C) \, d\Omega = -\int_\Omega \nabla \psi \cdot D(\mathbf{x}) \, \nabla C \, d\Omega + \int_{\Gamma} \psi (D(\mathbf{x}) \, \nabla C \cdot \mathbf{n}) \, d\Gamma
    \label{Eq:IntByParts}
\end{equation}

\noindent Step(4): Note that the flux is $\mathbf{J} = D(\mathbf{x}) \nabla C$. Some parts of the boundary will have known flux $\mathbf{J}_{\Gamma}$ (Neumann boundary conditions), while some other parts will have prescribed concentration $C_{\Gamma}$ (Dirichlet boundary conditions). Thus, the surface integral of Eq.~(\ref{Eq:IntByParts}) becomes 

\begin{equation}
    \int_{\Gamma} \psi (\mathbf{J} \cdot \mathbf{n}) \, d\Gamma = \int_{\Gamma} \psi (\mathbf{J}_{\Gamma} \cdot \mathbf{n}) \, d\Gamma 
\end{equation}

\noindent The weak form of Eq.~(\ref{Eq:Example}) becomes

\begin{equation}
    \begin{split}
    \int_\Omega \psi (\frac{\partial C}{\partial t}) -\int_\Omega \nabla \psi \cdot D(\mathbf{x}) \nabla C \, d\Omega + \int_{\Gamma} \psi (\mathbf{J}_{\Gamma} \cdot \mathbf{n}) \, d\Gamma = 0
    \end{split}
\end{equation}

\section{Heat transfer}

\subsection{Constitutive relations}

The transient heat transfer is expressed as

\begin{equation}
    \rho c \, \frac{\partial T}{\partial t} - \nabla \cdot (\mathbf{k} \, \nabla T) - Q = 0
\end{equation}

\noindent Where $\rho$ is the mass density, $c$ is the heat capacity, $Q$ is the heat generated per unit volume per unit time\footnote{source/sink term} and $\mathbf{k}$ is the conductivity matrix

\begin{equation}
    \mathbf{k} = \left[
        \begin{array}{ccc}
        k_x & 0 & 0 \\
        0 & k_y & 0 \\
        0 & 0 & k_z \\
        \end{array}\right]
\end{equation}

\subsection{The weak form}

The weak form, using $\mathbf{J} = -\mathbf{k} \nabla T$, is expressed as

\begin{equation}
    \int_\Omega \delta T \, \rho c \, \frac{\partial T}{\partial t} \, d \Omega + \int_\Omega \nabla \delta T \, \mathbf{k} \, \nabla T  \, d\Omega + \int_{\Gamma} \delta T \, (\mathbf{J}_{\Gamma} \cdot \mathbf{n}) \, d\Gamma - \int_\Omega \delta T \, Q  \,d\Omega = 0
    \label{Eq:HeatWeakForm}
\end{equation}

\subsection{FE discretization}

Using 

\begin{equation}
    T(\boldsymbol{\xi}) = \sum_i^n N_i(\boldsymbol{\xi})T_i \qquad 
    \nabla T(\boldsymbol{\xi}) = \sum_i^n \mathbf{B}_i(\boldsymbol{\xi})T_i
\end{equation}

\noindent Substituting in Eq.~(\ref{Eq:HeatWeakForm})\footnote{To simplify the notation, $\sum_i^n N_i(\boldsymbol{\xi}) \,T_i\equiv \mathbf{N} \mathbf{T}$ is used with the understanding that $\mathbf{T}$ is the vector of element nodal values.},  and eliminating $\delta T$

\begin{equation}
    \int_\Omega \mathbf{N}^\top \rho c \, \mathbf{N} \frac{\partial \mathbf{T}}{\partial t} \, d \Omega + \int_\Omega \mathbf{B}^\top \mathbf{k} \, \mathbf{B} \, \mathbf{T} \, d\Omega + \int_{\Gamma} \mathbf{N}^\top (\mathbf{J}_{\Gamma} \cdot \mathbf{n}) \, d\Gamma - \int_\Omega \mathbf{N}^\top Q \, d\Omega = 0
    \label{Eq:HeatDiscritized}
\end{equation}

\noindent Where the capacity matrix $\mathbf{M}$

\begin{equation}
    \mathbf{M} = \int_\Omega \mathbf{N}^\top \rho c \, \mathbf{N} \, d\Omega
\end{equation}

\noindent And the conductivity matrix $\mathbf{K}_\mathrm{D}$

\begin{equation}
    \mathbf{K}_\mathrm{D} = \int_\Omega \mathbf{B}^\top \mathbf{k} \, \mathbf{B} \, d\Omega
\end{equation}

\noindent The RHS load vector $\mathbf{F}$ for a constant source term $Q$

\begin{equation}
    \mathbf{F} = - \int_{\Gamma} \mathbf{N}^\top (\mathbf{J}_{\Gamma} \cdot \mathbf{n}) \, d\Gamma + \int_\Omega \mathbf{N}^\top Q \, d\Omega
\end{equation}

\noindent Then, Eq.~(\ref{Eq:HeatDiscritized}) can then be expressed as

\begin{equation}
    \mathbf{M} \, \dot{\mathbf{T}} + \mathbf{K}_\mathrm{D} \mathbf{T} = \mathbf{F}
    \label{Eq:TransientHeat}
\end{equation}

\noindent The implicit time integration of Eq.~(\ref{Eq:TransientHeat}) can be done as 

\begin{equation}
    (\mathbf{M} + \Delta t \, \mathbf{K}_\mathrm{D})\mathbf{T}^{t+1} = \Delta t \, \mathbf{F} + \mathbf{M} \, \mathbf{T}^t
\end{equation}

\section{Fully kinetic models for solute transport in metals}
\label{Sec: SoluteTransport}

The term \emph{fully kinetic} is used to describe this new class of models for hydrogen (and solute species in general) transport, as no local equilibrium assumptions are invoked in the formulation. Compared to the source-sink terms in the reaction-kinetics framework, accumulation or segregation in the fully kinetic framework is driven by diffusive mechanism. The core idea is to use a spatial field to identify microstructure features, which is used in the formulation to model interactions with these features. The order-parameters of the phase-field method, or functions of these order parameters, are chosen because they are bounded with well-defined maximum and minimum values and vary smoothly between them, which offers both mathematical and computational convenience \cite{Hussein2024, Hussein2024b, Hussein2025a}. The multiphase-field method of Stienbach et al. \cite{Steinbach2009} is chosen for running the phase-field based representative volume elements (RVEs). In the case of plastically deforming solids, the normalized dislocation density is used as a function of the equivalent plastic strain \cite{Hussein2025b}. For clarity and readability, the constitutive relations for segregation at microstructure features will be presented separately from that for dislocations, as the latter is typically coupled with phase-field fracture for damage evolution, which requires a special treatment and slight modification to the implementation. 

\subsection{Constitutive relations for interfaces and secondary phases segregation}

The total flux for hydrogen transport including interaction at the interfaces $g(\phi)_\mathrm{intf}$ and phases $\phi_n$ \cite{Hussein2024, Hussein2024b, Hussein2025a}

\begin{equation}
    \mathbf{J} = - \mathbf{D} \Big(\nabla c - \frac{\zeta_\mathrm{intf} \, c }{RT} \nabla g(\phi)_\mathrm{intf} - \frac{\zeta_n \, c }{RT} \nabla \phi_n \Big)
\end{equation}

\noindent An example algorithm for calculating the interface function $g(\phi)_\mathrm{intf}$ for a two-phase microstructure $\phi_i$ and $\phi_j$ is shown in Box~(\ref{Box: DoubleObstacle})\footnote{An implementation of this algorithm can be found in the supporting code \href{https://github.com/ahcomat/PHIMATS_RVE}{PHIMATS\_RVE} in \texttt{CaseStudies/Two\_Phase}.}

\nextboxnumber
\begin{boxH}[label=Box: DoubleObstacle]{Algorithm for calculating the double-obstacle potential for the three different interfaces resulting from a two-phase microstructure.}
    \begin{algorithm*}[H]
    \textbf{Initialize:} $g(\phi)_{ii} \gets 0$, $g(\phi)_{ij} \gets 0$,  $g(\phi)_{jj} \gets 0$ \;
    Define sets of grains: $\phi_i$, $\phi_j$ \;
    
    \For{$\alpha \in \{1, \dots, N\}$}{
        \For{$\beta \in \{1, \dots, i-1\}$}{

            \textbf{1.} Compute the interaction term:
            \[
            \text{term} \gets \phi_\alpha \cdot \phi_\beta
            \]
            
            \textbf{2.} Update $g(\phi)$ contributions based on grain categories:
            \begin{itemize}
                \item If $\alpha, \beta \in \phi_i$:
                \[
                g(\phi)_{ii} \gets g(\phi)_{ii} + \text{term}
                \]
                \item Else if $\alpha, \beta \in \phi_j$:
                \[
                g(\phi)_{jj} \gets g(\phi)_{jj} + \text{term}
                \]
                \item Else if $\alpha \in \phi_i$ and $\beta \in \phi_j$ (or vice versa):
                \[
                g(\phi)_{ij} \gets g(\phi)_{ij} + \text{term}
                \]
            \end{itemize}
        }
    }
    \textbf{Return:} Updated $g(\phi)$ terms: $g(\phi)_{ii}, g(\phi)_{ij}, g(\phi)_{jj}$
\end{algorithm*}
\end{boxH}

\noindent Applying the mass conservation 

\begin{equation}
    \frac{\partial c}{\partial t} - \nabla \cdot \mathbf{D} \Big( \nabla \, c - \frac{\zeta_\mathrm{intf} \, c }{RT} \nabla g(\phi)_\mathrm{intf} - \frac{\zeta_n c }{RT} \nabla \phi_n \Big) = 0
    \label{Eq:HydrogenDiffTrap}
\end{equation}

\subsubsection{The weak form}

Multiply Eq.~(\ref{Eq:HydrogenDiffTrap}) by a test function $\eta$ and integrate over the domain 

\begin{equation}
    \int_\Omega \eta \, \frac{\partial c}{\partial t} \, d\Omega - \int_\Omega \eta \, \nabla \cdot \mathbf{D} \Big( \nabla c - \frac{\zeta_\mathrm{intf} \, c }{RT} \nabla g(\phi)_\mathrm{intf} - \frac{\zeta_n \, c }{RT} \nabla \phi_n \Big) d\Omega = 0
\end{equation}

\noindent Integrate by parts the second term of the LHS and combining all terms, the weak form of Eq.~(\ref{Eq:HydrogenDiffTrap}) becomes

\begin{equation}
    \int_\Omega \eta \frac{\partial c}{\partial t}  d\Omega + \int_\Omega \nabla \eta \cdot \mathbf{D} \left( \nabla c - \frac{\zeta_\mathrm{intf} \, c }{RT} \nabla g(\phi)_\mathrm{intf} - \frac{\zeta_n \, c }{RT} \nabla \phi_n  \right)  d\Omega = - \int_{\Gamma} \eta \, (\mathbf{J}  \cdot \mathbf{n}) \, d\Gamma
    \label{Eq:TrapWeakForm}
\end{equation}

\subsubsection{FE discretization}

Using the integration-point values interpolated from the nodal points

\begin{equation}
    \begin{split}
        c = \mathbf{N \,c} \qquad 
        \nabla c = \mathbf{B} \, \mathbf{c}& \\
        \nabla g(\phi)_\mathrm{intf} = \mathbf{B} \, \boldsymbol{g(\phi)}_\mathrm{intf}& \qquad 
        \nabla \phi_n = \mathbf{B} \, \boldsymbol{\phi_n}
    \end{split}
\end{equation}

\noindent Substituting in \ref{Eq:TrapWeakForm} and eliminating $\eta$

\begin{equation}
    \begin{split}
        \int_\Omega \mathbf{N}^\top  \frac{\partial}{\partial t} \mathbf{N}\,\mathbf{c} \, d\Omega 
        &+ \int_\Omega \mathbf{B}^\top \mathbf{D} \mathbf{B} \, \mathbf{c} \, d\Omega
        - \int_\Omega \mathbf{B}^\top  \mathbf{D} \frac{\zeta_\mathrm{intf}}{RT} \mathbf{B} \, \boldsymbol{g(\phi)}_\mathrm{intf} \, \mathbf{N} \, \mathbf{c} \, d\Omega \\
        &- \int_\Omega \mathbf{B}^\top \mathbf{D} \frac{\zeta_{n}}{RT} \mathbf{B} \, \boldsymbol{\phi_{n}} \, \mathbf{N} \, \mathbf{c} \, d\Omega 
        = - \int_{\Gamma} \mathbf{N}^\top (\mathbf{J} \cdot \mathbf{n}) \, d\Gamma
        \end{split}
    \label{Eq:TrapDiscritized}
\end{equation}

\noindent Where the mass matrix $\mathbf{M}$

\begin{equation}
    \mathbf{M} = \int_\Omega \mathbf{N}^\top \mathbf{N} \, d\Omega
\end{equation}

\noindent And the diffusivity matrix $\mathbf{K}_\mathrm{D}$

\begin{equation}
    \mathbf{K}_\mathrm{D} = \int_\Omega \mathbf{B}^\top \mathbf{D} \mathbf{B} d\Omega
\end{equation}

\noindent And the interaction matrix $\mathbf{K}_\mathrm{I}$

\begin{equation}
    \mathbf{K}_\mathrm{I} = \int_\Omega \mathbf{B}^\top \mathbf{D}\frac{\zeta_\mathrm{intf}}{RT} \mathbf{B} \, \boldsymbol{g(\phi)}_\mathrm{intf} \mathbf{N} \, d\Omega + \int_\Omega \mathbf{B}^\top \mathbf{D}\frac{\zeta_{n}}{RT} \mathbf{B} \, \boldsymbol{\phi_{n}} \mathbf{N} \,d\Omega
\end{equation}

\noindent Eq.~(\ref{Eq:TrapDiscritized}) can then be expressed as

\begin{equation}
    \mathbf{M} \dot{\mathbf{c}} + [\mathbf{K}_\mathrm{D} - \mathbf{K}_\mathrm{I}] \mathbf{c} = \mathbf{F}
    \label{Eq:TrapTransient}
\end{equation}

\noindent The implicit time integration of Eq.~(\ref{Eq:TrapTransient}) can be done as 

\begin{equation}
    \frac{\mathbf{c}^{t+1}-\mathbf{c}^{t}}{\Delta t} \mathbf{M} + \mathbf{K}_\mathrm{D} \mathbf{c}^{t+1} - \mathbf{K}_\mathrm{I} \mathbf{c}^{t+1} = \mathbf{F}
\end{equation}

\begin{equation}
    [\mathbf{M} + \Delta t \, \mathbf{K}_\mathrm{D} - \Delta t \, \mathbf{K}_\mathrm{I}]\mathbf{c}^{t+1} = \Delta t \, \mathbf{F} + \mathbf{M} \mathbf{c}^t
    \label{Eq:TrapFEM}
\end{equation}

\noindent Note that Eq.~(\ref{Eq:TrapFEM}) is a linear system of equations and can be solved with either direct or iterative solvers. 

\subsection{Constitutive relations for interactions with dislocations and damage}

The hydrogen flux under mechanical loading is expressed as \cite{Hussein2025b}

\begin{equation}
  \mathbf{J} = -  \mathbf{D} \left(\nabla c - \frac{\bar{V}_\mathrm{H} \, c}{RT} \nabla \sigma_\mathrm{h} - \frac{\zeta_\rho \, c}{RT} \nabla \bar{\rho} \right)
  \label{Eq: ConFlux}
\end{equation}

\noindent where $\sigma_\mathrm{h}$ is the hydrostatic stress, $\nabla \bar{\rho}$ is the normalized dislocation density. Note that these values are calculated at the integration points. To calculate their gradients in PHIMATS, the integration point values are mapped to the nodes by volume averaging the contributions from all the integration points of the elements surrounding these nodes. Since these values are calculated as post-processing after the solution step, a staggered solution scheme is the preferred choice in PHIMATS when coupled with mechanics and phase-field fracture solvers. The mass conservation can be expressed as 

\begin{equation}
    s \frac{\partial c}{\partial t} - \nabla \cdot \mathbf{D} \, \Big(\nabla c - \frac{\bar{V}_\mathrm{H} \, c}{RT} \nabla \sigma_\mathrm{h} - \frac{\zeta_\rho \, c}{RT} \nabla \bar{\rho} \Big) + \frac{Z_d}{t} \, \phi^2 \, \Big(c - c_\mathrm{eq}\Big) = 0
    \label{Eq: MassConDis}
\end{equation}

\noindent where $s$ is the segregation capacity. In the case of damage, represented by the phase-field variable $\phi$, the crack surface acts as a source/sink for hydrogen, which is represented by the last term of Eq.~(\ref{Eq: ConFlux}). $Z_\mathrm{d}$ is the source/sink strength, $t$ is the time and $c_\mathrm{eq}$ is the equilibrium concentration calculated from

\begin{equation}
    c_\mathrm{eq} = c_B \exp\Big( \frac{\bar{V}_\mathrm{H} \, \sigma_\mathrm{h} + \zeta_\rho \, \bar{\rho}}{RT} \Big)
    \label{Eq: Conb}
\end{equation}

\noindent where $c_B$ is a constant boundary concentration value. 

\subsubsection{The weak form}

Multiply Eq.~(\ref{Eq: MassConDis}) by a test function $\eta$ and integrate over the domain 

\begin{equation}
  \begin{split}
  \int_\Omega s\, \eta \, \frac{\partial c}{\partial t} \, d\Omega 
  + \int_\Omega \nabla \eta \cdot \mathbf{D} \Big( 
      \nabla c - \frac{\bar{V}_\mathrm{H} \, c}{RT} \nabla \sigma_\mathrm{h} 
       & - \frac{\zeta_\rho \, c}{RT} \nabla \bar{\rho} \Big) \, d\Omega 
  + \int_\Omega \frac{Z_d}{t} \, \phi^2 \, \eta \, c \, d\Omega = \\
  & - \int_{\Gamma} \eta \, \mathbf{J} \cdot \mathbf{n} \, d\Gamma 
  + \int_\Omega \frac{Z_d}{t} \, \phi^2 \, \eta \, c_\mathrm{eq} \, d\Omega
  \end{split}
  \label{Eq: MassCon}
\end{equation}

\subsubsection{FE discretization}

Using the integration-point values interpolated from the nodal points

\begin{equation}      
  c = \mathbf{N \,c}  \qquad 
  \nabla c = \mathbf{B} \, \mathbf{c} \qquad 
  \nabla \sigma_\mathrm{h} = \mathbf{B} \, \boldsymbol{\sigma}_\mathrm{h} \qquad
  \nabla \bar{\rho} = \mathbf{B} \, \bar{\boldsymbol{\rho}} \qquad
  \label{Eq: DerivMat}
\end{equation}

\noindent Substituting in \ref{Eq: MassCon} and eliminating $\eta$

\begin{equation}
  \begin{aligned}
    \int_\Omega s\, \mathbf{N}^\top\mathbf{N} \frac{\partial\mathbf{c}}{\partial t} \, d\Omega 
    &+ \int_\Omega \mathbf{B}^\top \mathbf{D} \mathbf{B} \, \mathbf{c} \, d\Omega \,
    - \int_\Omega \mathbf{B}^\top \mathbf{D} \frac{\bar{V}_\mathrm{H}}{RT} \mathbf{B} \, \boldsymbol{\sigma}_\mathrm{h} \, \mathbf{N} \, \mathbf{c} \, d\Omega \,-\int_\Omega \mathbf{B}^\top \mathbf{D} \frac{\zeta_\rho}{RT} \mathbf{B} \, \bar{\boldsymbol{\rho}} \, \mathbf{N} \, \mathbf{c} \, d\Omega  \\
    & +\int_\Omega \mathbf{N}^\top\mathbf{N} \frac{Z_d}{t} \, \phi^2 \mathbf{c} \, d\Omega = - \int_{\Gamma} \mathbf{N}^\top \mathbf{J} \cdot \mathbf{n} \, d\Gamma \, +\int_\Omega \mathbf{N}^\top\mathbf{N} \frac{Z_d}{t} \, \phi^2 \mathbf{c}_\mathrm{eq} \, d\Omega
  \end{aligned}
  \label{Eq: Discret}
\end{equation}

\noindent Defining the matrices

\begin{equation}
  \begin{split}
    &\mathbf{M} = \int_\Omega s\, \mathbf{N}^\top \mathbf{N} \, d\Omega \qquad
    \mathbf{K_\mathrm}_\mathrm{d} = \int_\Omega \mathbf{B}^\top \mathbf{D} \mathbf{B} \, d\Omega \\
    &\mathbf{K}_\mathrm{I} = \int_\Omega \mathbf{B}^\top \mathbf{D}\frac{\bar{V}_\mathrm{H}}{RT} \mathbf{B} \, \boldsymbol{\sigma}_\mathrm{h} \mathbf{N} \, d\Omega + \int_\Omega \mathbf{B}^\top \mathbf{D}\frac{\zeta_\rho}{RT} \mathbf{B} \, \bar{\boldsymbol{\rho}} \mathbf{N} \,d\Omega \\
    & \mathbf{K}_\mathrm{S} = \int_\Omega \mathbf{N}^\top\mathbf{N} \frac{Z_d}{t} \, \phi^2 \, d\Omega
  \end{split}
\end{equation}

\noindent And the nodal RHS load vector 

\begin{equation}
    \mathbf{F} = - \int_{\Gamma} \mathbf{N}^\top \mathbf{J} \cdot \mathbf{n} \, d\Gamma  \, +\int_\Omega \mathbf{N}^\top\mathbf{N} \frac{Z_d}{t} \, \phi^2 \mathbf{c}_\mathrm{eq} \, d\Omega
\end{equation}

\noindent Substituting in Eq. (\ref{Eq: Discret})

\begin{equation}
    \mathbf{M} \, \dot{\mathbf{c}} + (\mathbf{K}_\mathrm{D} - \mathbf{K}_\mathrm{I} + \mathbf{K}_\mathrm{S}) \, \mathbf{c} = \mathbf{F}
    \label{Eq: TimeDisc}
\end{equation}

\noindent Using the implicit Euler time integration to Eq. (\ref{Eq: TimeDisc}), the FE linear system of equations becomes

\begin{equation}
    [\mathbf{M} + \Delta t \, (\mathbf{K}_\mathrm{D} - \, \mathbf{K}_\mathrm{I} + \mathbf{K}_\mathrm{S})] \, \mathbf{c}^{t+1} = \Delta t \, \mathbf{F} + \mathbf{M} \, \mathbf{c}^t
    \label{Eq: TrapTransient}
\end{equation}

\section{Isotropic linear elasticity -- Small strain formulation}

\subsection{Constitutive relations}

The stress-strain relations for linear elasticity are

\begin{equation}
    \begin{split}
    \boldsymbol{\sigma} &= \mathbf{C} \boldsymbol{\varepsilon} \\
    \sigma_{ij} &= C_{ijkl} \varepsilon_{kl}
    \end{split}
\end{equation}

\noindent For isotropic linear elasticity, the stress-strin relations are expressed using the Lam\'e constants $\lambda = E\nu/(1+\nu)$ and $\mu = G = E/2(1+\nu)$

\begin{equation}
    \boldsymbol{\sigma} = 2G \boldsymbol{\varepsilon} + \lambda \Tr(\boldsymbol{\varepsilon}) \mathbf{I}
\end{equation}

\noindent And in component form

\begin{equation}
    \left[
        \begin{array}{rrr}
        \sigma_{xx} &  \sigma_{xy} & \sigma_{xz} \\
        \sigma_{yx}  &  \sigma_{yy}  &  \sigma_{yz} \\
        \sigma_{zx} &  \sigma_{zy} & \sigma_{zz} \\
        \end{array}
    \right] = 2G \left[
        \begin{array}{rrr}
        \varepsilon_{xx} &  \varepsilon_{xy} & \varepsilon_{xz} \\
        \varepsilon_{yx}  &  \varepsilon_{yy}  &  \varepsilon_{yz} \\
        \varepsilon_{zx} &  \varepsilon_{zy} & \varepsilon_{zz} \\
        \end{array}
    \right] + \lambda(\varepsilon_{xx}+\varepsilon_{yy}+\varepsilon_{zz}) \left[
        \begin{array}{rrr}
        1 &  0 & 0 \\
        0 &  1 & 0 \\
        0 &  0 & 1 \\
        \end{array}
    \right]
\end{equation}

\noindent For numerical convenience, it is more practical to represent the stress-strain relations in \emph{Voigt} notation, where the shear components are stored as \emph{engineering} shears.

\begin{equation}
    \begin{split}
    \left[
        \begin{array}{r}
        \sigma_{xx} \\
        \sigma_{yy} \\
        \sigma_{zz} \\
        \sigma_{xy} \\
        \sigma_{yz} \\
        \sigma_{xz} \\
        \end{array}
    \right] &= \left[
        \begin{array}{cccccc}
        \lambda+2\mu & \lambda & \lambda & 0 & 0 & 0 \\
        \lambda & \lambda+2\mu & \lambda & 0 & 0 & 0 \\
        \lambda & \lambda & \lambda+2\mu & 0 & 0 & 0 \\
        0 & 0 & 0 & \mu & 0 & 0 \\
        0 & 0 & 0 & 0 & \mu & 0 \\
        0 & 0 & 0 & 0 & 0 & \mu \\
        \end{array}
    \right] \left[
        \begin{array}{c}
        \varepsilon_{xx} \\
        \varepsilon_{yy} \\
        \varepsilon_{zz} \\
        2\varepsilon_{xy} \\
        2\varepsilon_{yz} \\
        2\varepsilon_{xz} \\
        \end{array}
    \right] \\
    &= \left[
        \begin{array}{cccccc}
        \lambda+2\mu & \lambda & \lambda & 0 & 0 & 0 \\
        \lambda & \lambda+2\mu & \lambda & 0 & 0 & 0 \\
        \lambda & \lambda & \lambda+2\mu & 0 & 0 & 0 \\
        0 & 0 & 0 & \mu & 0 & 0 \\
        0 & 0 & 0 & 0 & \mu & 0 \\
        0 & 0 & 0 & 0 & 0 & \mu \\
        \end{array}
    \right] \left[
        \begin{array}{c}
        \varepsilon_{xx} \\
        \varepsilon_{yy} \\
        \varepsilon_{zz} \\
        \gamma_{xy} \\
        \gamma_{yz} \\
        \gamma_{xz} \\
        \end{array}
    \right] \\
    \end{split}
\end{equation}

\noindent For plane-strain

\begin{equation}
    \left[
        \begin{array}{r}
        \sigma_{xx} \\
        \sigma_{yy} \\
        \sigma_{xy} \\
        \end{array}
    \right] = \frac{E}{(1+\nu)(1-2\nu)}\left[
        \begin{array}{ccc}
        1-\nu & \nu & 0 \\
        \nu & 1-\nu & 0 \\
        0 & 0 & \frac{1-2\nu}{2} \\
        \end{array}
    \right] \left[
        \begin{array}{c}
        \varepsilon_{xx} \\
        \varepsilon_{yy} \\
        \gamma_{xy} \\
        \end{array}
    \right] \\
\end{equation}

\noindent $\sigma_{zz}\neq 0$\footnote{In PHIMATS, $\sigma_{zz}$ is an output for plasticity models. However, it is not an output for elasticity models, yet, it could be evaluated in post-processing.} and is calculated from

\begin{equation}
    \sigma_{zz} = \nu(\sigma_{xx} + \sigma_{yy})
\end{equation}

\noindent And for plane-stress

\begin{equation}
    \left[
        \begin{array}{r}
        \sigma_{xx} \\
        \sigma_{yy} \\
        \sigma_{xy} \\
        \end{array}
    \right] = \frac{E}{(1-\nu^2)}\left[
        \begin{array}{ccc}
        1 & \nu & 0 \\
        \nu & 1 & 0 \\
        0 & 0 & \frac{1-\nu}{2} \\
        \end{array}
    \right] \left[
        \begin{array}{c}
        \varepsilon_{xx} \\
        \varepsilon_{yy} \\
        \gamma_{xy} \\
        \end{array}
    \right] \\
\end{equation}

\noindent Note that $\varepsilon_{zz}\neq 0$ and is calculated from

\begin{equation}
    \varepsilon_{zz} = -\frac{1}{3E}(\sigma_{xx} + \sigma_{yy})
\end{equation}

\noindent Accordingly, the hydrostatic stress $\sigma_\mathrm{h}$ for plane-strain 

\begin{equation}
    \sigma_{\mathrm{h}} = \frac{1 + \nu}{3} (\sigma_{xx} + \sigma_{yy})
\end{equation}

\noindent And for plane-stress

\begin{equation}
    \sigma_{\mathrm{h}} = \frac{1}{3} (\sigma_{xx} + \sigma_{yy})
\end{equation}

\subsection{The weak form}

The quasi-static balance of linear momentum is 

\begin{equation}
    \begin{split}
    \nabla \cdot \boldsymbol{\sigma} + \mathbf{b} &= 0 \\
    \frac{\partial \sigma_{ij}}{\partial x_j} + b_i &= 0 \\
    \end{split}
    \label{Eq:LinearMomentum}
\end{equation}

\noindent The problem statement for quasi-static linear elasticity is find $u_i$, $\varepsilon_{ij}$ and $\sigma_{ij}$ that satisfies

\begin{enumerate}
    \item The strain-displacement relation $\varepsilon_{ij} = 1/2(\partial u_i/\partial x_j + \partial u_j/\partial x_i)$
    \item The stress-strain relation $\sigma_{ij} = C_{ijkl}\varepsilon_{kl}$
    \item The quasi-static linear balance of momentum $\frac{\partial \sigma_{ij}}{\partial x_j} + b_i = 0$
    \item The displacement and traction boundary conditions $u_i = u_i^{\Gamma_1}$ on $\Gamma_1$ and $\sigma_{ij}n_j = t_j^{\Gamma_2}$ on $\Gamma$
\end{enumerate}

\noindent The principle of virtual work is an equivalent integral form of the linear balance of momentum PDE. This makes it more convenient for numerical solution using computers. This essentially involves all the 4-steps mentioned in Section \ref{Sec: WeakForm}. We multiply the strong form Eq.~(\ref{Eq:LinearMomentum}) by a kinematically admissible virtual displacement field $\delta u_i$ (test function), which means it satisfies $\delta u_i=0$ on $\Gamma^1$. Following Bower \cite{Bower2009}, \emph{This is a complicated way for saying that the small perturbation displacement $\delta u_i$ satisfies the boundary conditions.} The associated \emph{virtual} strain field 

\begin{equation}
    \delta \varepsilon_{ij} = 1/2(\partial u_i/\partial x_j + \partial u_j/\partial x_i)
\end{equation}

\noindent The weak form is then

\begin{equation}
    \int_\Omega \delta \varepsilon_{ji} \sigma_{ij}  d\Omega - \int_\Omega \delta u_j b_i  d\Omega - \int_{\Gamma} \delta u_j t_i  d\Gamma = 0
    \label{Eq:LinElasWeakForm}
\end{equation}

\begin{equation}
    \int_\Omega \delta \boldsymbol{\varepsilon}^{\top} \boldsymbol{\sigma} d\Omega - \int_\Omega \delta \mathbf{u}^{\top} \mathbf{b}  d\Omega - \int_{\Gamma} \delta \mathbf{u}^{\top} \mathbf{t}  d\Gamma = 0
    \label{Eq:LinElasWeakForm}
\end{equation}

\noindent Will satisfy the strong form for all possible $\delta \mathbf{u}$.
    
\subsection{FE discretization}

The displacements $\mathbf{u}(\boldsymbol{\xi})$ within any point in the element $\boldsymbol{\xi}$ is obtained from the nodal values of the element using the isoparametric formalism

\begin{equation}
    \mathbf{u}(\boldsymbol{\xi}) = \left[
        \begin{array}{c}
        u_{x}(\boldsymbol{\xi}) \\
        u_{y}(\boldsymbol{\xi}) \\
        u_{z}(\boldsymbol{\xi}) \\
        \end{array}
    \right] = \sum_i^n N_i(\boldsymbol{\xi})u_i
    \label{Eq:IsoparamDisp}
\end{equation}



\noindent The strains can be obtained from the derivatives of the nodal displacements. For simplicity, we use a 2D case. The small-strain displacement in Voigt notation becomes

\begin{equation}
    \boldsymbol{\varepsilon} = \left[
        \begin{array}{c}
        \varepsilon_{xx} \\
        \varepsilon_{yy} \\
        \gamma_{xy} \\
        \end{array}
        \right] = \left[\begin{array}{c}
        \frac{\partial u_x}{\partial x} \\
        \frac{\partial u_y}{\partial y} \\
        \frac{\partial u_x}{\partial y} + \frac{\partial u_y}{\partial x} \\
        \end{array}\right]
\end{equation}

\noindent Using the displacement values from Eq.~(\ref{Eq:IsoparamDisp}), the strains could be evaluated using the cartesian derivatives of the shape functions similar to Eq.~(\ref{Eq:IsoparamDeriv}). However, due to the special case of the Voigt representation of the strain tensor, the cartesian derivatives of the shape functions are collected in a matrix called the $\mathbf{B}$ matrix\footnote{Note that $\mathbf{B}$ is of size $\texttt{nStres} \times \texttt{nElDispDofs}$, where $\texttt{nStre}$ is the number of stress components and $\texttt{nElDispDofs}$ is the number of displacement degrees of freedom per element. This will have to be evaluated for every gauss point $\boldsymbol{\xi}$, i.e. stored in size \texttt{nGaus}.}, which has a special form as will be described later. Therefore 

\begin{equation}
    \boldsymbol{\varepsilon}(\boldsymbol{\xi}) = \sum_i^n \mathbf{B}_i(\boldsymbol{\xi}) u_i
\end{equation}

\nextboxnumber
\begin{boxH}[label=Bmatrix]{Example for calculating strain in 4-node quad element}
    \begin{equation*}
        \begin{split}
            \left[
                \begin{array}{c}
                \varepsilon_{xx} \\
                \varepsilon_{yy} \\
                \gamma_{xy} \\
                \end{array}
            \right] = \left[\begin{array}{cccccccc}
                \frac{\partial N_1}{\partial x} & 0 & \frac{\partial N_2}{\partial x} & 0 & \frac{\partial N_3}{\partial x} & 0 & \frac{\partial N_4}{\partial x} & 0 \\
                0 & \frac{\partial N_1}{\partial y} & 0 & \frac{\partial N_2}{\partial y} & 0 & \frac{\partial N_3}{\partial y} & 0 & \frac{\partial N_4}{\partial y} \\
                \frac{\partial N_1}{\partial y} & \frac{\partial N_1}{\partial x} & \frac{\partial N_2}{\partial y} & \frac{\partial N_2}{\partial x} & \frac{\partial N_3}{\partial y} & \frac{\partial N_3}{\partial x} & \frac{\partial N_4}{\partial y} & \frac{\partial N_4}{\partial x} \\
                \end{array}\right]  \left[\begin{array}{c}
                    u1_{x} \\
                    u1_{y} \\
                    u2_{x} \\
                    u2_{y} \\
                    u3_{x} \\
                    u3_{y} \\
                    u4_{x} \\
                    u4_{y} \\
                    \end{array}
                    \right]
        \end{split}
    \end{equation*}
\end{boxH}

\noindent Using

\begin{equation}
    \delta\mathbf{u} = \sum_i^n N_i(\boldsymbol{\xi})\delta \mathbf{u}
\end{equation}

\noindent And 

\begin{equation}
    \delta \boldsymbol{\varepsilon} = \sum_i^n \mathbf{B}_i(\boldsymbol{\xi}) \delta\mathbf{u}
\end{equation}

\noindent The weak form in Eq.~(\ref{Eq:LinElasWeakForm}) can be written as

\begin{equation}
    \int_\Omega \mathbf{B}^\top \delta \mathbf{u}^{\top} \boldsymbol{\sigma}  d\Omega - \int_\Omega \mathbf{N}^\top \delta \mathbf{u}^{\top} \mathbf{b} d\Omega - \int_{\Gamma} \mathbf{N}^\top \mathbf{u}^{\top} \mathbf{t} d\Gamma = 0
\end{equation}

\noindent Using the stress-strain relations and eliminating the arbitrary $\delta \mathbf{u}$ 

\begin{equation}
    \int_\Omega \mathbf{B}^\top \boldsymbol{\sigma} d\Omega - \int_\Omega \mathbf{N}^\top \mathbf{b} d\Omega - \int_{\Gamma} \mathbf{N}^\top \mathbf{t} d\Gamma = 0
\end{equation}

\begin{equation}
    \int_\Omega \mathbf{B}^\top \mathbf{C} \, \mathbf{B} \mathbf{u} d\Omega - \int_\Omega \mathbf{N}^\top \mathbf{b} d\Omega - \int_{\Gamma} \mathbf{N}^\top \mathbf{t} d\Gamma = 0
\end{equation}

\noindent Ignoring body forces

\begin{equation}
    \int_\Omega \mathbf{B}^\top \mathbf{C} \, \mathbf{B} \mathbf{u} d\Omega - \int_{\Gamma} \mathbf{N}^\top \mathbf{t} d\Gamma = 0
\end{equation}

\noindent Or

\begin{equation}
    \mathbf{f}_\mathrm{int} - \mathbf{f}_\mathrm{ext} = 0
    \label{Eq: StressEquilibrium}
\end{equation}

\noindent The stiffness matrix $\mathbf{K}$

\begin{equation}
    \mathbf{K} = \int_\Omega \mathbf{B}^\top \mathbf{C} \, \mathbf{B} d\Omega
\end{equation}

\noindent Since the domain is discretized in \texttt{nElem} elements and introducing the element stiffness matrix $\mathbf{k}^{(e)}$, the \emph{global} stiffness matrix becomes

\begin{equation}
    \mathbf{K} = \sum^\texttt{nElem} \mathbf{k}^{(e)}
\end{equation}

\noindent The element stiffness matrix can be defined as 

\begin{equation}
    \mathbf{k}^{(e)} = \int_{\Omega_{(e)}}{\Big[\sum_i^n \mathbf{B}_i(\boldsymbol{\xi})\Big]}^\top \mathbf{C} \sum_i^n \mathbf{B}_i(\boldsymbol{\xi}) d\Omega_{(e)}
\end{equation}

\noindent Where the integral can be evaluated using Gauss quadrature. The reaction force can then be calculated from 

\begin{equation}
    \mathbf{f}_\mathrm{reaction} = -\mathbf{f}_\mathrm{int} = -\int_\Omega \mathbf{B}^\top \boldsymbol{\sigma} d\Omega 
\end{equation}

\section{Non-linear finite element method}

\subsection{Tangential stiffness matrix}

The nonlinear (material or geometric) behavior will result in a nonlinear system of equations Eq.~(\ref{Eq: StressEquilibrium}) which can be solved using the Newton-Raphson method. This requires the linearization of the stress-strain relation via \cite{Borst2012}

\begin{equation}
    \delta \boldsymbol{\sigma} = \frac{\partial \boldsymbol{\sigma}}{\partial \boldsymbol{\varepsilon}} \delta \boldsymbol{\varepsilon} = \mathbf{C} \, \delta \boldsymbol{\varepsilon}
\end{equation}

\noindent Where $\mathbf{C}$ is the tangential stiffness matrix.

\subsection{The Newton-Raphson method for nonlinear systems}

The Newton-Raphson method is used for solving nonlinear equations. Its implementation is based on the evaluation of the nonlinear function and its derivative. The Taylor expansion of $f(x + \Delta x)$ around $x$ is 

\begin{equation}
    f(x + \Delta x) = f(x) + f' \Delta x + \cdots
\end{equation}

\noindent Where $f'=df(x)/dx$. The root of the function could be reached by the iteration scheme:

\begin{equation}
    x_{n+1} = x_n - \frac{f(x_n)}{f'(x_n)}
\end{equation}

\subsection{Incremental iterative analysis}

The balance of momentum equation of load increment can be written as 

\begin{equation}
    \begin{split}
    \mathbf{f}^{(n+1)}_\mathrm{ext} - \mathbf{f}^{(n)}_\mathrm{int} - \Delta \mathbf{f}_\mathrm{int} &= 0 \\
    \mathbf{f}^{(n+1)}_\mathrm{ext} - \mathbf{f}^{(n)}_\mathrm{int} - \mathbf{K} \Delta \mathbf{u} &= 0 \\
    \mathbf{f}^{(n+1)}_\mathrm{ext} - \mathbf{f}^{(n)}_\mathrm{int} - \mathbf{R} &= 0 \\
    \end{split}
\end{equation}

\noindent The system becomes nonlinear due to the nonlinear material behavior, which is solved using the Newton-Raphson method. The process for each loading increment is summarized in Algorithm (\ref{Alg: NRSmallPlast}) \cite{Borst2012}

\nextboxnumber
\begin{boxH}[label=Box: DoubleObstacle]{Newton-Raphson algorithm for a nonlinear system of equations}
\begin{algorithm*}[H]
    \textbf{Initialize:} $j \gets 0$, $\mathbf{u}^{(0)}$ \;
    \While{$\| \mathbf{R}^{(j)} \| > \text{Tolerance}$ \textbf{and} $j < \text{Max Iterations}$}{
        \textbf{1.} Calculate the total strain: 
        \[
        \boldsymbol{\varepsilon}^{(j)} = \sum \mathbf{B} \mathbf{u}^{(j)} 
        \]
        
        \textbf{3.} Calculate the internal force: 
        \[
        \mathbf{f}_\mathrm{int}^{(j)} = \int_\Omega \mathbf{B}^\top \boldsymbol{\sigma}^{(j)} \, d\Omega
        \]

        \textbf{4.} Calculate the residual:
        \[
        \mathbf{R}^{(j)} = \mathbf{f}_\mathrm{ext} - \mathbf{f}_\mathrm{int}^{(j)}
        \]

        \textbf{5.} Assemble the tangent stiffness matrix:
        \[
        \mathbf{K}^{(j)} = \int_\Omega \mathbf{B}^\top \mathbf{C} \, \mathbf{B} \, d\Omega
        \]

        \textbf{6.} Solve the linear system:
        \[
        \Delta \mathbf{u}^{(j+1)} = \mathbf{K}^{(j)^{-1}} \mathbf{R}^{(j)}
        \]

        \textbf{7.} Update displacement:
        \[
        \mathbf{u}^{(j+1)} = \mathbf{u}^{(j)} + \Delta \mathbf{u}^{(j+1)}
        \]

        \textbf{8.} Check Convergence:
        \If{$\| \mathbf{R}^{(j+1)} \| \leq \text{Tolerance}$}{
            \textbf{Break}
        }
    }
    \textbf{Return:} Converged displacement $\mathbf{u}^{(j+1)}$
    \label{Alg: NRSmallPlast}
\end{algorithm*}
\end{boxH}

\section{$J_2$ Plasticity -- Isotropic hardening}

\subsection{Constitutive relations}

\subsubsection{Von Mises stress}

The hydrostatic stress is defined as 

\begin{equation}
    \sigma_\mathrm{h} = \sigma_{ii} = \Tr(\boldsymbol{\sigma}) =  (\sigma_{xx}+\sigma_{yy}+\sigma_{zz})
\end{equation}

\noindent The deviatoric stress is defines as

\begin{equation}
    \begin{split}
        \boldsymbol{\sigma}' &= \sigma_{ij}-\sigma_{ii}\delta_{ij} = \boldsymbol{\sigma} - \Tr(\boldsymbol{\sigma})\mathbf{I} \\
        &= \left[
        \begin{array}{ccc}
        \sigma_{xx}' &  \sigma_{xy} & \sigma_{xz} \\
        \sigma_{yx}  &  \sigma_{yy}'  &  \sigma_{yz} \\
        \sigma_{zx} &  \sigma_{zy} & \sigma_{zz}' \\
        \end{array}
    \right] = \left[
        \begin{array}{ccc}
        \sigma_{xx}-\sigma_\mathrm{h} &  \sigma_{xy} & \sigma_{xz} \\
        \sigma_{yx}  &  \sigma_{yy}-\sigma_\mathrm{h}  &  \sigma_{yz} \\
        \sigma_{zx} &  \sigma_{zy} & \sigma_{zz}-\sigma_\mathrm{h} \\
        \end{array}
    \right]
\end{split}
\end{equation}

\noindent The von Mises, or, \emph{equivalent} stress is written as 

\begin{equation}
    \begin{split}
    \sigma_\mathrm{eq} &= \sqrt{\frac{1}{2}\Big[(\sigma_{xx}-\sigma_{xx})^2+(\sigma_{yy}-\sigma_{zz})^2+(\sigma_{xx}-\sigma_{zz})^2\Big]+3(\sigma_{xy}^2+\sigma_{yz}^2+\sigma_{xz}^2)} \\
    &= \sqrt{\sigma_{xx}^2+\sigma_{yy}^2+\sigma_{zz}^2-\sigma_{xx}\sigma_{yy}-\sigma_{yy}\sigma_{zz}-\sigma_{zz}\sigma_{xx}+3(\sigma_{xy}^2+\sigma_{yz}^2+\sigma_{xz}^2)} \\
    &= \sqrt{\frac{3}{2}\sigma_{ij}\sigma_{ij}-\frac{1}{2}\sigma_{kk}^2} \\
    &= \sqrt{\frac{3}{2}\sigma_{ij}'\sigma_{ij}'} \\
    \end{split}
\end{equation}

\subsubsection{Plastic incompressibility}

\noindent The total strain is the sum of the elastic and plastic strains as 

\begin{equation}
    \begin{split}    
    \boldsymbol{\varepsilon} &= \boldsymbol{\varepsilon}^\mathrm{e} + \boldsymbol{\varepsilon}^\mathrm{p} \\
    \varepsilon_{ij} &= \varepsilon_{ij}^\mathrm{e} + \varepsilon_{ij}^\mathrm{p}
    \end{split}
\end{equation}

\noindent Since there is no volume change during plastic deformation, the consequence for the plastic strain rate

\begin{equation}
    \dot{\varepsilon}_{xx}^\mathrm{p}+\dot{\varepsilon}_{yy}^\mathrm{p}+\dot{\varepsilon}_{zz}^\mathrm{p} = 0
\end{equation}

\noindent The equivalent plastic strain rate is defined as

\begin{equation}
    \dot{p} =\sqrt{\frac{2}{3}\dot{\boldsymbol{\varepsilon}}^\mathrm{p}:\dot{\boldsymbol{\varepsilon}}^\mathrm{p}} =  \sqrt{\frac{2}{3}\dot{\varepsilon}^\mathrm{p}_{ij}\dot{\varepsilon}^\mathrm{p}_{ij}}
\end{equation}

\noindent It should be noted that the plastic strain rates are deviatoric by construction due to the incompressibility condition. The accumulated plastic strain is defined as 

\begin{equation}
    p = \int_0^t = \dot{p} \, dt
\end{equation}

\subsubsection{Yield criterion}

The von Mises yield function $f(\boldsymbol{\sigma}, \varepsilon^\mathrm{eq})$ defines a yield surface $f(\boldsymbol{\sigma}, \varepsilon^\mathrm{eq})=0$

\begin{equation}
    \begin{split}
    f &= \sigma_\mathrm{eq} - \sigma_\mathrm{y} - H(p) < 0: \quad \text{Elastic deformation} \\
    f &= \sigma_\mathrm{eq} - \sigma_\mathrm{y} - H(p) = 0: \quad \text{Plastic deformation} \\
    \end{split}
\end{equation}

\noindent Where $H(p)$ is the strain hardening function.

\subsubsection{Hardening laws implemented in PHIMATS}

\begin{itemize}
    \item Power law hardening \texttt{[PowerLaw]}. 
    \item Kocks-Mecking-Estrin \texttt{[KME]}.
\end{itemize}

\subsubsection{Flow rule}

The flow rule determines the direction of plastic flow \cite{Dunne2005, Doghri2000}. This is accomplished through the \emph{normality condition}. It states that, in the principal stress space, the increment of plastic strain tensor $d\boldsymbol{\varepsilon}^\mathrm{p}$ is normal to the tangent of the yield surface at the loading point $\frac{\partial f}{\partial \boldsymbol{\sigma}}$, thus 

\begin{equation}
    \dot{\boldsymbol{\varepsilon}}^\mathrm{p} = \dot{\lambda} \frac{\partial f}{\partial \boldsymbol{\sigma}} \qquad d\boldsymbol{\varepsilon}^\mathrm{p} = d\lambda \frac{\partial f}{\partial \boldsymbol{\sigma}}  
\end{equation}

\noindent Where $\dot{\lambda}$ is the plastic multiplier, which determines the magnitude of plastic strain rate. For a von Mises material, the plastic multiplier is the increment of equivalent plastic strain 

\begin{equation}
    \dot{p} = \dot{\lambda} \qquad dp = d\lambda
\end{equation}

\noindent Thus, for a von Mises material, the normality condition can be stated as 

\begin{equation}
    d\boldsymbol{\varepsilon}^\mathrm{p} = d\lambda \frac{\partial f}{\partial \boldsymbol{\sigma}} = \frac{3}{2} \, dp \frac{\boldsymbol{\sigma'}}{\sigma_\mathrm{eq}}
\end{equation}

\noindent With the tensor normal to the yield surface

\begin{equation}
    \boldsymbol{N} = \frac{3}{2} \frac{\boldsymbol{\sigma'}}{\sigma_\mathrm{eq}} \quad N_{ij} = \frac{3}{2} \frac{\sigma'_{ij}}{\sigma_\mathrm{eq}}
\end{equation}

\subsubsection{Consistency condition}

The consistency condition states that during plastic deformation, the load point, or the solution, should stay on the yield surface. This condition enables determining the plastic multiplier, or for a von Mises plasticity, the equivalent plastic strain. Mathematically, it could be stated as   

\begin{equation}
    f = \sigma_\mathrm{eq} - \sigma_\mathrm{y} - H(p) = 0
\end{equation}

\subsection{Return mapping algorithm}

The total elastic strain at increment $(n+1)$ is written as 

\begin{equation}
    \boldsymbol{\varepsilon}^\mathrm{e}_{(n+1)} = (\boldsymbol{\varepsilon}_{(n)} + \Delta \boldsymbol{\varepsilon}) - (\boldsymbol{\varepsilon}^\mathrm{p}_{(n)} + \Delta \boldsymbol{\varepsilon}^\mathrm{p})
\end{equation}

\noindent The stress is then

\begin{equation}
    \boldsymbol{\sigma}_{(n+1)} = 2G (\boldsymbol{\varepsilon}_{(n)} + \Delta \boldsymbol{\varepsilon} - \boldsymbol{\varepsilon}^\mathrm{p}_{(n)}-\Delta \boldsymbol{\varepsilon}^\mathrm{p}) + \frac{1}{3}\lambda \Tr(\boldsymbol{\varepsilon}_{(n)} + \Delta \boldsymbol{\varepsilon} - \boldsymbol{\varepsilon}^\mathrm{p}_{(n)}-\Delta \boldsymbol{\varepsilon}^\mathrm{p}) \mathbf{I}
\end{equation}

\noindent And since $\Tr(\Delta \boldsymbol{\varepsilon}^\mathrm{p}) = 0$

\begin{equation}
    \boldsymbol{\sigma}_{(n+1)} = \underbrace{
        2G (\boldsymbol{\varepsilon}_{(n)} + \Delta \boldsymbol{\varepsilon}- \boldsymbol{\varepsilon}^\mathrm{p}_{(n)}) + \frac{1}{3}\lambda \Tr(\boldsymbol{\varepsilon}_{(n)} + \Delta \boldsymbol{\varepsilon}- \boldsymbol{\varepsilon}^\mathrm{p}_{(n)}) \mathbf{I}}_{\substack{\text{elastic predictor} \\ \text{(trial stress $\boldsymbol{\sigma}^\mathrm{tr}$)}}}
        \underbrace{
        -2G \Delta \boldsymbol{\varepsilon}^\mathrm{p}}_{\text{plastic corrector}}
        \label{Eq: Predictor}
\end{equation}

\noindent Eq.~(\ref{Eq: Predictor}) could be written as 

\begin{equation}
    \begin{split}    
    \boldsymbol{\sigma}_{(n+1)} = \boldsymbol{\sigma}^\mathrm{tr} - 2G \Delta \boldsymbol{\varepsilon}^\mathrm{p} 
    &= \boldsymbol{\sigma}^\mathrm{tr} - 2G \Delta p \boldsymbol{N} \\
    &= \boldsymbol{\sigma}^\mathrm{tr} - 2G \Delta p \frac{3}{2}\frac{{\boldsymbol{\sigma}'}^\mathrm{tr}}{\sigma_\mathrm{eq}}
    \end{split}
\end{equation}

\noindent Since only the deviatoric part of the stress needs to be computed. Omitting the subscript $_{(n+1)}$ for simplicity

\begin{equation}
    \begin{split}
    \boldsymbol{\sigma}' &= {\boldsymbol{\sigma}'}^\mathrm{tr} - 2G \Delta p \boldsymbol{N} \\
    \frac{2}{3} {\sigma_\mathrm{eq}}\boldsymbol{N} &= \frac{2}{3}{\sigma_\mathrm{eq}}^\mathrm{tr}\boldsymbol{N} - 2G \Delta p \boldsymbol{N} \\
    {\sigma_\mathrm{eq}} &= \sigma_\mathrm{eq}^\mathrm{tr} - 3G \Delta p
    \end{split}
\end{equation}

\noindent Using the yield function

\begin{equation}
    f = \sigma_\mathrm{eq}^\mathrm{tr} - 3G \Delta p \underbrace{- H(p) - \sigma_\mathrm{y}}_{\text{$\sigma_\mathrm{eq}$}} = 0
\end{equation}

\noindent Therefore, the problem reduces to solving a nonlinear scalar equation for the unknown $\Delta p$. Using Taylor's expansion 

\begin{equation}
    f + \frac{df}{d\Delta p} d\Delta p + \cdots = 0 , \, \text{with} \, \frac{df}{d\Delta p}= - 3G - \frac{dH(p)}{dp} 
\end{equation}

\noindent We can write the increment of the equivalent plastic strain as 

\begin{equation}
    d\Delta p = \frac{\sigma_\mathrm{eq}^\mathrm{tr} - 3G \Delta p - H(p) - \sigma_\mathrm{y}}{3G + \frac{dH(p)}{dp}}
\end{equation}

\noindent Now, we can solve for the iterations $(k)$

\begin{equation}
    \begin{split}
    d\Delta p &= \frac{f^{(k)}}{3G + \frac{dH(p^{(k)})}{dp}}  \\ 
    \Delta p^{(k+1)} & = \Delta p^{(k)} + d\Delta p \\  
    p^{(k+1)} &= p^{(k)} + \Delta p^{(k+1)} \\
    f^{(k+1)} &= \sigma_\mathrm{eq}^\mathrm{tr} - 3G \Delta p - H(p) - \sigma_\mathrm{y}
    \end{split}
    \label{Eq: RM}
\end{equation}

\noindent Once a convergence criterion $|f^{(k+1)}\le\mathrm{tol}|$ is achieved, the solution at $(n+1)$ can be evaluated as

\begin{equation}
    \begin{split}
    \boldsymbol{\varepsilon}^\mathrm{p} &= \frac{3}{2} \, \Delta p \boldsymbol{N} = \frac{3}{2} \, \Delta p \frac{{\boldsymbol{\sigma}'}^\mathrm{tr}}{\sigma_\mathrm{eq}}  \\
    \boldsymbol{\varepsilon}^\mathrm{e}&= \boldsymbol{\varepsilon} - \boldsymbol{\varepsilon}^\mathrm{p} \\
    \boldsymbol{\sigma} &= 2G \boldsymbol{\varepsilon}^\mathrm{e} + \lambda \Tr(\boldsymbol{\varepsilon}^\mathrm{e}) \mathbf{I} 
    \end{split}
\end{equation}

\subsection{Tangential stiffness matrix}

For isotropic hardening plasticity, $\mathbf{\mathbf{C}_\mathrm{ep}}$ can be calculated from 

\begin{equation}
    \mathbf{C}_\mathrm{ep} = \mathbf{C}_\mathrm{e} - \frac{\mathbf{C}_\mathrm{e} : \boldsymbol{N} \otimes \boldsymbol{N} : \mathbf{C}_\mathrm{e}}{\frac{2}{3} H + \boldsymbol{N} : \mathbf{C}_\mathrm{e} : \boldsymbol{N}}
    \label{Eq: IsoHardTangStiffness}
\end{equation}

\section{Phase-field fracture}

\subsection{Introduction to the geometric phase-field fracture}

\begin{figure}[hbt!]
    \includegraphics[width=8cm]{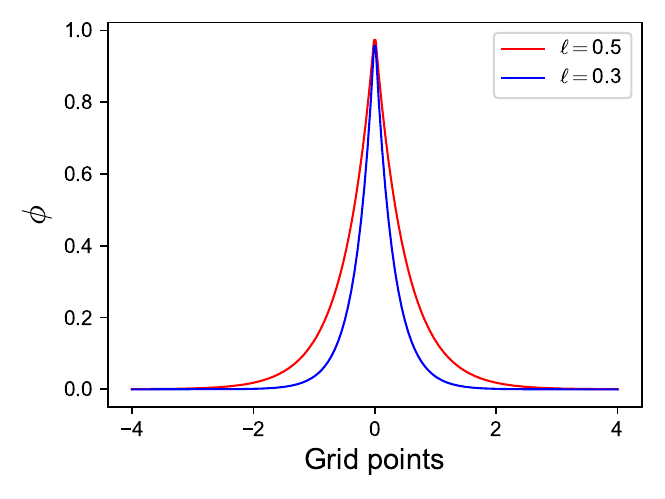}
    \centering
    \caption{Phase-field representation of the 1D diffuse crack topology.} 
    \label{Fig: 1DPF}
  \end{figure}

The idea of phase-field fracture is similar to the continuum (diffuse) damage mechanics, where a sharp crack is approximated with the phase-field function \cite{Biner2017}. The diffuse (smeared) crack topology can be represented in 1D by the phase-field order parameter $\phi$ as shown in Fig.~(\ref{Fig: 1DPF})

\begin{equation}
    \phi(x) = \exp\Big(-\frac{|x|}{\ell}\Big)
\end{equation}

\noindent where $\ell$ is a parameter determining the spreading (smearing) of the crack, with $\ell\rightarrow0$ is the sharp interface limit. $\phi$ varies between 0 for intact material and 1 for fully damaged material. The geometric phase-field fracture of Miehe et al \cite{Miehe2015_p1, Miehe2015_p2, Miehe2016}, represents the regularized crack surface density functional as

\begin{equation}
    \Gamma_\ell(\phi) = \int_V \, \gamma(\phi, \nabla \phi) dV, \qquad \text{with} \qquad \gamma_\ell(\phi, \nabla \phi) = \frac{1}{2\ell} (\phi^2 + \ell^2 |\nabla \phi|^2).
    \label{Eq: CrackSurfDen}
\end{equation}

\noindent The crack evolution is governed by the balance of a crack driving force $\mathcal{H}$ and the crack geometric resistance 

\begin{equation}
    \phi - \ell^2 \Delta \phi = (1 - \phi) \, \mathcal{H}
    \label{Eq: CrackBalance}
\end{equation}

\noindent The crack driving force $\mathcal{H}$ is expressed as a function of a \emph{defined state} and ensures irreversibility, such that

\begin{equation}
    \mathcal{H} = \max_{s \in [0, \, t]} \tilde{\mathcal{D}}(state, s)
\end{equation}

\noindent The advantage of the geometric phase-field fracture formulation is that it allows flexibility and modularity in defining $\tilde{\mathcal{D}}(state)$, i.e., the crack driving force. In PHIMATS, the implemented crack driving forces are defined below.

\subsection{Definition of crack driving force $\tilde{\mathcal{D}}(state)$ options available in PHIMATS}

\subsubsection*{Brittle}

\begin{equation}
    \tilde{\mathcal{D}}^\mathrm{b} = \frac{\tilde{\psi}^\mathrm{e+}}{w_\mathrm{c}}
\end{equation}

\noindent where $w_\mathrm{c}$ is the critical work energy density, $\tilde{\psi}^\mathrm{e+}$ is the positive part of the strain energy density, defined as \cite{Miehe2010_IJNME, Miehe2010_CMAME}

\begin{equation}
    \tilde{\psi}^{\mathrm{e}+}
    = \frac{1}{2}\lambda \langle \mathrm{tr}(\boldsymbol{\varepsilon}^\mathrm{e}) \rangle^2_+ + \mu \mathrm{tr}([\boldsymbol{\varepsilon}^\mathrm{e+}]^2) \quad  \text{and} \quad 
    \tilde{\psi}^{\mathrm{e}-}
    = \frac{1}{2}\lambda \langle \mathrm{tr}(\boldsymbol{\varepsilon}^\mathrm{e}) \rangle^2_- + \mu \mathrm{tr}([\boldsymbol{\varepsilon}^\mathrm{e-}]^2) 
\end{equation}

\noindent which are calculated from the spectral decomposition of the elastic strain 

\begin{equation}
    \boldsymbol{\varepsilon}^{\mathrm{e}+}
    = \sum_{i=1}^{n} \langle \varepsilon^\mathrm{e}_i \rangle_+ \, \mathbf{n}_i \otimes \mathbf{n}_i \, \, ,
    \qquad \text{and} \qquad
    \boldsymbol{\varepsilon}^{\mathrm{e}-}
    = \sum_{i=1}^{n} \langle \varepsilon^\mathrm{e}_i \rangle_- \, \mathbf{n}_i \otimes \mathbf{n}_i \, .
\end{equation}

\noindent In a PHIMATS driver code, the \emph{spectral decomposition} and the \emph{brittle} crack driving force are called from the member functions of a \texttt{PFFModel} object

\begin{equation*}
    \begin{split}
    \texttt{PFFModel::CalcPsiSpectral(} &\texttt{vector<BaseElemPFF*> pffElem,} \\
    &\texttt{vector<BaseElemMech*> mechElem,} \\
    &\texttt{vector<BaseMechanics*> mechMat)}
    \end{split}
\end{equation*}

\begin{equation*}
    \texttt{PFFModel::CalcDrivFrocElas(vector<BaseElemPFF*> pffElem)}
\end{equation*}

\subsubsection*{Elastoplastic}

\begin{equation}
    \tilde{\mathcal{D}}^\mathrm{ep} = \frac{\tilde{\psi}^\mathrm{e+} + \tilde{w}^\mathrm{p}}{w_\mathrm{c}} 
\end{equation}

\noindent where $\tilde{w}^\mathrm{p}$ is the plastic work density defined as 

\begin{equation}
    \tilde{w}^\mathrm{p} = \int_0^t \boldsymbol{\tilde{\sigma}} \!:\! \dot{\boldsymbol{\varepsilon}}^\mathrm{p} \, ds
\end{equation}

\noindent where $\boldsymbol{\tilde{\sigma}}$ is the \emph{undegraded} stress. The \emph{elastoplastic} crack driving force is called from

\begin{equation*}
    \begin{split}
    \texttt{PFFModel::CalcDrivFrocEP(}&\texttt{vector<BaseElemPFF*> pffElem,} \\
    &\texttt{vector<BaseElemMech*> mechElem)}
    \end{split}
\end{equation*}

\subsubsection*{Elastoplastic with a threshold}

\begin{equation}
    \langle \tilde{\mathcal{D}}^\mathrm{ep} \rangle = \zeta \left\langle \frac{\tilde{\psi}^\mathrm{e+} + \tilde{w}^\mathrm{p}}{\tilde{w}_\mathrm{c}} - 1 \right\rangle
\end{equation}

\noindent where $\zeta$ is a parameter that controls the post-critical behavior \cite{Miehe2015_p1}. It can be seen that this form ensures that crack initiation cannot develop under the condition $\tilde{\psi}^\mathrm{e+} + \tilde{w}^\mathrm{p}<\tilde{w}_\mathrm{c}$, i.e., the elastoplastic crack driving force is less than the critical value. The elastoplastic crack driving force with a \emph{threshold} is called from

\begin{equation*}
    \begin{split}
    \texttt{PFFModel::CalcDrivFrocEP\_TH(} &\texttt{vector<BaseElemPFF*> pffElem,} \\ 
    &\texttt{vector<BaseElemMech*> mechElem,} \\
    &\texttt{const double zeta)}
    \end{split}
\end{equation*}

\noindent Note that \texttt{zeta} is an optional function argument with a default value of 1.

\subsection{Stress degradation}

The degradation accompanied by crack evolution is described by 

\begin{equation}
  g(\phi)_\mathrm{d} = (1-\phi)^2
\end{equation}

\noindent The evolution of stress follows the conservation of linear momentum 

\begin{equation}
  \nabla \boldsymbol{\sigma} = 0, \qquad \boldsymbol{\sigma} = g(\phi)_\mathrm{d} \, \tilde{\boldsymbol{\sigma}} = g(\phi)_\mathrm{d} \, \mathbf{C}_\mathrm{e} \!:\! \boldsymbol{\varepsilon}^\mathrm{e} \, .
\end{equation}

\subsection{Modifications to the return mapping algorithm during damage}

In the case of a degrading elastoplastic material, the return mapping algorithm described in Section \ref{Eq: RM} is modified to account for the damage according to

\begin{equation}
    \begin{split}
    d\Delta p &= \frac{f^{(k)}}{3G \, g(\phi)_\mathrm{d} + \frac{dH(p^{(k)})}{dp}}  \\ 
    \Delta p^{(k+1)} & = \Delta p^{(k)} + d\Delta p \\  
    p^{(k+1)} &= p^{(k)} + \Delta p^{(k+1)} \\
    f^{(k+1)} &= \sigma_\mathrm{eq}^\mathrm{tr} - 3G \, g(\phi)_\mathrm{d} \, \Delta p - H(p) - \sigma_\mathrm{y}
    \end{split}
    \label{Eq: RMPFF}
\end{equation}

\noindent Together with the degraded elastoplastic tangent matrix 

\begin{equation}
  \mathbf{C}_\mathrm{ep} = g(\phi)_\mathrm{d} \, \tilde{\mathbf{C}}_\mathrm{ep} 
\end{equation}

\subsection{The weak form and FEM discretization}

The weak form of Eq.~(\ref{Eq: CrackBalance}) by a test function $\eta_\phi$

\begin{equation}
  \int_V \phi \, \eta_\phi - \ell^2 \nabla \phi \cdot \nabla \eta_\phi \, dV - \int_V (1 - \phi) \mathcal{H} \, \eta_\phi \, dV = 0
\end{equation}

\noindent Using the shape and derivative functions

\begin{equation}      
  \phi = \mathbf{N \,\phi} \, ,  \qquad 
  \nabla \phi = \mathbf{B} \, \mathbf{\phi} \, . 
  \label{Eq: DerivMatPFF}
\end{equation}

\noindent The residual can be obtained by eliminating the arbitrary test function 

\begin{equation}
    \mathbf{R}_\phi =  \int_V (\mathbf{N}\boldsymbol{\phi} \, \mathbf{N}^\top  + \ell^2 \, \mathbf{B}^\top \mathbf{B} \, \boldsymbol{\phi}) dV - \int_V \mathcal{H} (1-\mathbf{N}\boldsymbol{\phi}) \, \mathbf{N}^\top dV
    \label{Eq: PhiResidual}
\end{equation}

\noindent By decoupling the mechanics and the crack evolution problems \cite{Miehe2010_CMAME}, a semi-implicit time integration scheme for updating $\mathcal{H}$ results in the linear system of equations \cite{Miehe2015_p1}

\begin{equation}
    \int_V \Big[ (\mathcal{H}+1)\mathbf{N}^\top \mathbf{N}  + \ell^2 \, \mathbf{B}^\top \mathbf{B} \Big] \, dV \boldsymbol{\phi} - \int_V \mathcal{H} \, \mathbf{N}^\top \, dV = 0
    \label{Eq: PFFDiscrete}
\end{equation}

\begin{equation}
    \mathbf{K}_{\phi} \boldsymbol{\phi} = \mathbf{f}_\mathcal{H} 
\end{equation}

\printbibliography

\end{document}